%% file: wl1.tex
\documentclass{article}
\usepackage{amsfonts}
\usepackage{amsmath}
\usepackage{graphicx,epsfig}
\usepackage{amssymb}
\usepackage{rotating}
\usepackage{theorem}
\newtheorem{prop}{Proposition}[section]
\theorembodyfont{\rmfamily}
\usepackage{float}
\usepackage[bf]{caption}

\begin{document}
\thispagestyle{plain}
\begin{center}
\begin{Large}
\textbf{Polyhedral Deformations of Cone Manifolds}\\
\vspace{2cm}
  A Aalam
\end{Large}
\end{center}
\vspace{1cm}

\textbf{Abstract} Two single parameter families of polyhedra $P(\psi)$ are
constructed in three dimensional spaces of constant curvature
$C(\psi)$. Identification of the faces of the polyhedra via isometries
results in cone manifolds $M(\psi)$ which are topologically $S^1\times
S^2$, $S^3$ or singular $S^2$ . The singular set of $M(\psi)$ can have self intersections for
some values of $\psi$ and can also be the Whitehead link or form
other configurations. Curvature varies continuously with $\psi$. At
$\psi=0$ spontaneous surgery occurs and the topological type of
$M(\psi)$ changes. This phenomenon is described.

\section*{0  Introduction}
We study continuous families of cone manifolds $M_\psi$ parametrised by cone
angle which begin at cone angle zero with the complement of the Whitehead link in
$S^3$. We consider the case of equal cone angles on all singular link
components. Increasing cone angles the families trace different
paths in Dehn surgery space joined by what we call a {\it Dehn surgery
transition point}.  The cone structures for certain non-zero
values of cone angles exist in projective models or in $S^3$. 
 
In one Dehn surgery
direction the cone manifolds are for certain cone angles, obtained by
surgery on the Whitehead link in $S^3$ resulting in a
topologically distinct singular set in $S^2 \times S^1$. As cone angle
is increased the topological type of the singular set changes and the
hyperbolic cone manifold develops two cusps and becomes $S^3$ at cone angle
$\frac{2}{3}$$\pi$. The topological type of the singular set and the
structure of $M_\psi$ remain unaltered as cone angle
increases beyond $\frac{2}{3}$$\pi$ until we reach a cone angle $\omega$
where $M_\psi$ becomes $\Bbb R^3$ with topologically
the same type of singularity. Increasing cone angle past $\omega$ the
singular set reverts back to its pre-$\frac{2}{3}$$\pi$ cone angle
topological type and $M_\psi$ becomes spherical in $S^2
\times S^1$. At cone angle $\pi$ the underlying polyhedron becomes a
lens in $S^3$ from which the cone manifold is obtained by suitable
identifications. For cone angles in the interval [$\pi, \zeta$], $M_\psi$ is spherical and the topological type of its singular set is
unchanged but it is now in $S^3$. At cone angle $\zeta$, $M_\psi$ becomes the
suspension of a sphere with four cone points. It remains the well
understood sphere with four singularities for cone angles larger than $\zeta$. 

Investigating the deformation on the other side of Dehn surgery we
obtain the Whitehead link in $S^3$ for certain non-zero cone
angles. A complete investigation will be carried out later.

Working in hyperbolic or spherical three space of constant curvature
it is usual for curvature to be normalised to plus or minus one. While
there are good reasons for this convention, this does not allow us to
envisage a continuous family of cone structures in which curvature
changes from positive to negative with Euclidean space a point in a
continuum.cf \cite{HLM}.

Here a {\it Cone Manifold} is a PL manifold with a possibly empty
codimension two locally flat submanifold called the {\it singular
  set}. In dimensions two and three the singular set consists of
isolated points and curves but not graphs respectively. The geometric
model is a spherical, Euclidean, or hyperbolic space of constant
curvature where the constant is any real number. Points in the
complement of the singular set have neighbourhoods homeomorphic to
neighbourhoods in the model. Points on the singular set have
neighbourhoods homeomorphic to neighbourhoods in the topological space
obtained by identifying boundaries of the intersection of two half
spaces referred to as a {\it wedge} in the model. The homeomorphism takes the singular set to the
axis of rotation in the topological space and transition functions are
isometries. {\it Orbifolds} are represented by discrete structures where cone
angles are of the form $\frac{2\pi}{n} , n \in \Bbb N$. cf \cite{HLM}

The topological space obtained by identifying in pairs faces 
 of a polyhedron in a space of constant curvature via isometries is a
 cone manifold if :

\begin{enumerate}
\item No edge is identified with its inverse in the equivalence class
induced by the identifications, and the identifications of wedges
along faces are cyclic for each equivalence class of edges.
\item The cone angle at each edge, i.e. the sum of the dihedral angles
  about the edge is $\leq 2\pi$.
\item The neighbourhood of each vertex is a cone on a sphere. If a
  vertex neighbourhood is, for example, a cone on a torus the
  topological space is not a manifold.
\item There are either two or no edges emanating from a vertex where
  cone angle is not 2$\pi$. If two edges then the cone angles must be
  equal and the edges must be lined up. The vertex neighbourhood
  is obtained by identifying wedge half planes.
\end{enumerate} 
cf \cite{HLM}.

It is reasonable to expect that our deformation process holds for
any hyperbolic link complement since the polyhedral description of a
link complement utilised here is canonical and we can
expect to be able to deform this construction by opening cusps in
$M_0$ in the way we have described. Our methods can be viewed as
part of an approach to describe all compact connected 3-manifolds through
deformations by removing singularities of branched singular covers of
$S^3$ along universal links thereby providing an approach to the Poincare Conjecture in dimension three.
The important phenomenon of self intersecting singular set in a cone
manifold occurs in our work. Thus our construction may prove useful in
studying this phenomenon.

\section{Cone angle interval $[0,\omega]$}\label{section:s1}

\begin{prop}
 $M_\psi$ is $S^2 \times S^1$ with two singular components and is
 hyperbolic for cone angles $\psi \in (0,\frac{2}{3}
 \pi)$. $M_{\frac{2}{3}\pi}$ is cusped hyperbolic and has a self intersecting
 singular set. There exists $\omega \in \Bbb R$ such that $M_\psi$ is
 hyperbolic with a self intersecting singular set when
 $\frac{2}{3}\pi < \psi < \omega$ and $M_\omega$ is $\Bbb R^3$.
 \end{prop}

{\bf Proof} $M_0$ the complement of the Whitehead link in $S^3$ is
obtained from a two cusped octahedron $P_0$ with identifications described
in the projective Klein model of $\Bbb H^3$ in figure
\ref{fig:fig3}. We refer the reader to \cite{Rat}, \cite{Thu2} and
\cite{Thu1} for details of this construction. 

\begin{figure}[H]
\includegraphics[height=35mm,width=40mm]{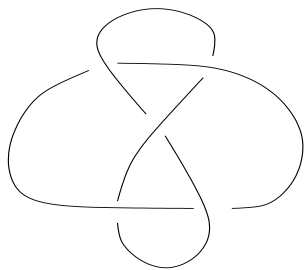}
\centering
\caption{The Whitehead link}\label{fig:fig2}
\end{figure}

\begin{figure}[H]
\includegraphics[height=70mm,width=70mm]{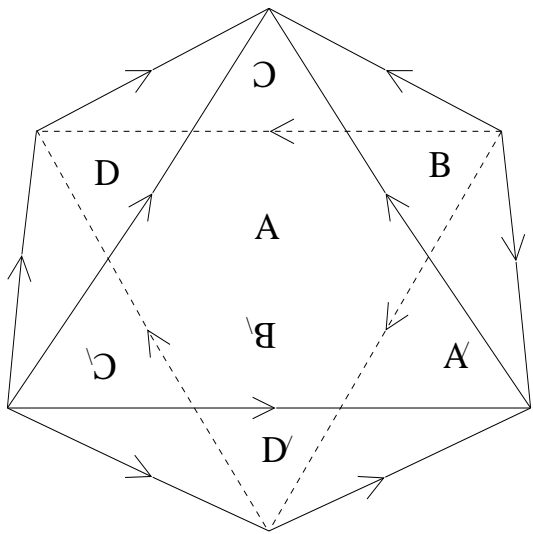}
\centering
\caption{$P_0$ with identifications}\label{fig:fig3}
\end{figure}

\begin{figure}[H]
\includegraphics[height=60mm]{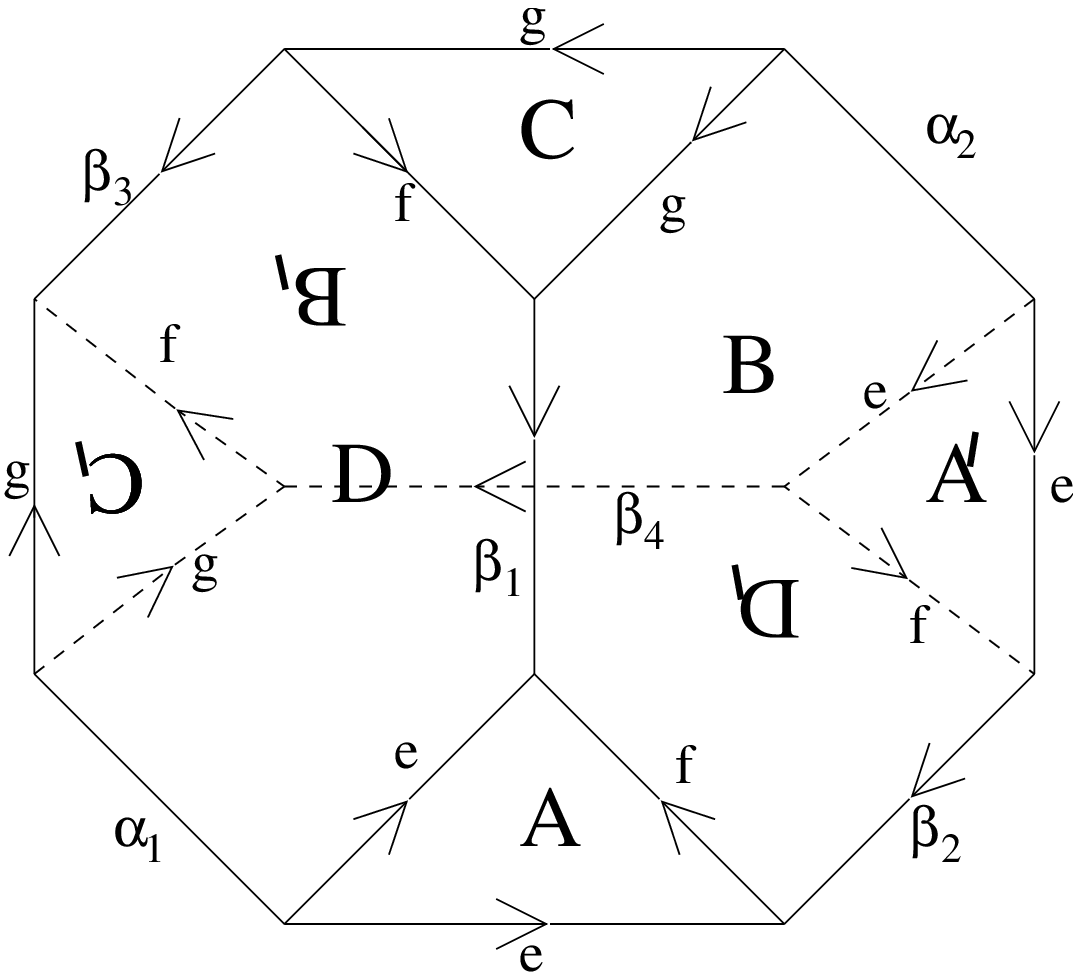}
\centering
\caption{opening cusps}\label{fig:fig4}
\end{figure}

\begin{figure}[H]
\includegraphics[height=80mm,width=80mm]{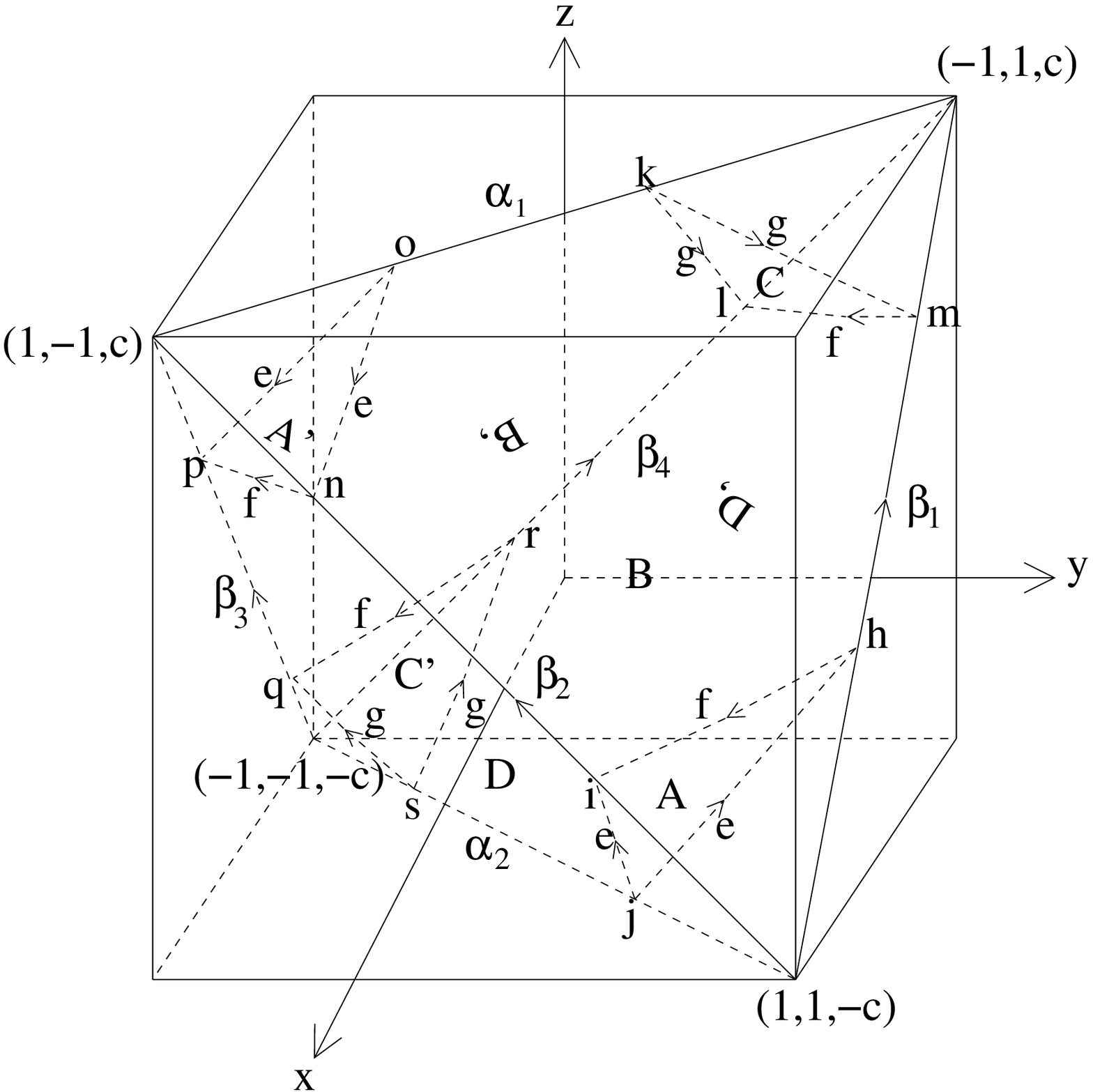}
\centering
\caption{opening cusps}\label{fig:fig5}
\end{figure}

To find out more about $P_\psi$ we view it relative to
a coordinate system and inside a {\it reference box} in the Klein
model as in figure \ref{fig:fig5}.  The length, width and
height of the reference box are denoted by $a, b$ and $c$ respectively. Opening cusps in $M_0$ we obtain $P_\psi$ for $\psi \in
(0,\frac{2}{3}\pi)$ as in figure \ref{fig:fig4}.

From figure \ref{fig:fig5} we have the following identifications:\\

Face identifications:
\begin{gather}
A \longleftrightarrow A' \hspace{5mm} i.e. \hspace{5mm} \bigtriangleup(onp)
\longleftrightarrow \bigtriangleup(jhi) \\ C \longleftrightarrow C'
\hspace{5mm} i.e. \hspace{5mm} \bigtriangleup(sqr) \longleftrightarrow
\bigtriangleup(klm)\\ D \longleftrightarrow D' \hspace{5mm}
i.e. \hspace{5mm} hexagon (inpqsj) \longleftrightarrow hexagon (hmlrsj)\\ B
\longleftrightarrow B'\hspace{5mm} i.e. \hspace{5mm} hexagon (nokmhi)
\longleftrightarrow hexagon (poklrq)
\end{gather} \\

Edge pairings:
\begin{gather}
ok \longleftrightarrow ok\label{5}\\  sj \longleftrightarrow sj\label{6}\\ on
\longleftrightarrow op \longleftrightarrow ji \longleftrightarrow jh\label{7}\\
sr\longleftrightarrow sq \longleftrightarrow km \longleftrightarrow
kl\label{8}\\ np \longleftrightarrow ml \longleftrightarrow hi
\longleftrightarrow rq\label{9}\\ in \longleftrightarrow hm \longleftrightarrow
lr \longleftrightarrow pq \label{10}
\end{gather}

Edge pairings \ref{5}, \ref{6} and \ref{10} represent singular
components of $M_\psi$. 

We observe: 
\begin{gather}
\notag\text{The dihedral angle between planes incident}\\
\text{in a member of (\ref{5}) or (\ref{6}) is the cone angle}\; \psi.\label{11}
\end{gather}
\\
Moreover (\ref{7}), (\ref{8}) and (\ref{9}) are not part of the singular
set. Therefore
\begin{gather}
\notag\text{The dihedral angle between planes incident}\\ \text{in a member of (\ref{7}), (\ref{8}) or (\ref{9}) is}\; \frac{\pi}{2}.\label{12}
\end{gather}

Information in (\ref{10}) represents segments $\beta_1, \beta_2,
\beta_3$ and $\beta_4$ which are identified to give a singular
component $\beta$ of the singular set of $M_\psi$. Considering
(\ref{10}) with equal cone angles on all singular components we
deduce:
\begin{gather}
\text{Planes incident in a member of (\ref{10}) intersect at dihedral
  angle}\; \frac{\psi}{4}.\label{13}
\end{gather} 

We deduce from (\ref{10}) that $a=b$. Box coordinates can therefore be
normalised to give $a=1=b$ and $c$ where $c$ is measured in the
$z-$direction. The reference box is therefore a cube with square base
and height $2c$ as in figure \ref{fig:fig5}.

\begin{figure}[H]
\includegraphics[height=70mm,width=70mm]{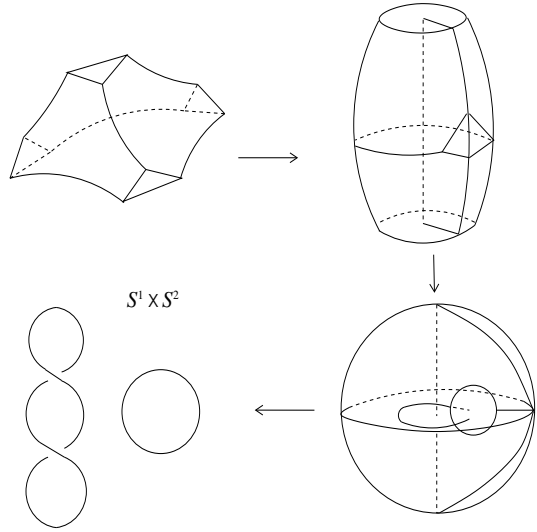}
\centering
\caption{Topology of $M_\psi$}\label{fig:OC}
\end{figure}

Topologically $M_\psi$ is $S^1 \times S^2$ with two unlinked singular
components as in figure \ref{fig:OC}.

To find out about the geometry of $M_\psi$ we look at $P_\psi$ inside
the Klein model of $\Bbb H^3$. Let $R$ be the radius of the Klein ball
$B_R$ and let $c$ represent the height of the main box. We will show
that $P_\psi$ lives inside $B_R$ and can be specified by $R$
and $c$. We also have
\begin{gather}
\text{curvature of}\; \Bbb H^3 = -\frac{1}{R^2}\label{curv}
\end{gather}

The geometry of $M_\psi$ can therefore be described in terms of
$\psi$.
We now obtain $R$ and $c$ as functions of $\psi$.

\begin{figure}[H]
\includegraphics[height=57mm,width=60mm]{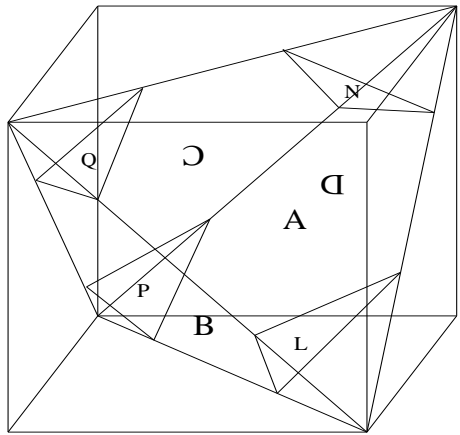}
\centering
\caption{$P_\psi$ in the Klein model}\label{fig:fig6}
\end{figure}

Figure \ref{fig:fig6} depicts $P_\psi$ in the Klein model. 

In homogeneous coordinates the equations of planes of $P_\psi$ and their poles are:

\begin{equation}\label{planes}
\begin{gathered}
plane \\ A \\ B \\ C \\ D
\end{gathered}
\begin{gathered}
equation \\ cx+cy+z=c \\ cx-cy-z=c \\ -cx-cy+z=c \\ -cx+cy-z=c
\end{gathered}\qquad
\begin{gathered}
pole \\
(cR,cR,R,c) \\ (cR,-cR,-R,c) \\ (-cR,-cR,R,c) \\ (-cR,cR,-R,c)
\end{gathered}
\end{equation}

Let ${\bf v}=(v_1,v_2,v_3,v_4), {\bf w}=(w_1,w_2,w_3,w_4) \in
\Bbb R^4$. Then
the hyperbolic bilinear form is given by 
\begin{equation}\label{HBF}
{\langle{\bf v},{\bf w} \rangle}_\Bbb H = v_1 w_1 + v_2 w_2 + v_3 w_3 - v_4 w_4
\end{equation}

Let $\theta$ be the dihedral angle of intersection between planes $P$
and $Q$ with poles ${\bf v}$ and ${\bf w}$ respectively. Then

\begin{equation}\label{angle}
\cos \theta = -\dfrac{{\langle {\bf v},{\bf w} \rangle}_\Bbb H}{\sqrt {{\langle
    {\bf v}, {\bf v} \rangle}_\Bbb H} \sqrt {{\langle {\bf w},{\bf w}
    \rangle}_\Bbb H}}
\end{equation}
 
Combining (\ref{10}), (\ref{13}), (\ref{planes}) and (\ref{angle}) we can
derive expressions for $\cos \psi$ and $\cos \dfrac{\psi}{4}$ in terms
of $R$ and $c$ to obtain :

\begin{gather}
c^2 = \dfrac{1+\cos\psi}{2\cos\dfrac{\psi}{4}-\cos\psi+1}\label{c2}\\ \notag
\\ R^2=\dfrac{1+\cos\psi}{2\cos\dfrac{\psi}{4}+\cos\psi-1}\label{R2}
\end{gather}

Substituting $(R, c)-$expressions of $\cos\dfrac{ \psi}{4}$ and
$\cos\psi$ in (\ref{c2}) and (\ref{R2}) we can verify these identities.

We note
$\omega=4\cos^{-1}(\dfrac{2}{\sqrt3}\cos(\dfrac{1}{3}\cos^{-1}(\dfrac{-3\sqrt3}{8})))\approx
2.311984...$ is the smallest positive value of $\psi$ for which the
denominator of (\ref{R2}) is zero. Hence $\omega$ is the smallest
positive value of $\psi$ for which $R^2$ is infinite. This combined
with (\ref{curv}) implies that curvature is zero and $M_\psi$ is
therefore Euclidean when $\psi=\omega$.

Combining (\ref{12}), (\ref{planes}), (\ref{HBF}) and (\ref{angle}) the equation of plane $N$ in figure \ref{fig:fig6} is :
\begin{equation}\label{N}
N : -x+y+cz=R^2
\end{equation}

Let ${\bf p}=(-1,1,-c,R)$ denote the pole of $N$ with $\Bbb H^3$
embedded in $\Bbb RP^3$. We have
\begin{equation}
\langle{\bf p},{\bf
  p}\rangle=2+c^2-R^2=\dfrac{3-4\cos^2\dfrac{\psi}{4}}{(4\cos^3\dfrac{\psi}{4}-4\cos\dfrac{\psi}{4}-1)(4\cos^3\dfrac{\psi}{4}-4\cos\dfrac{\psi}{4}+1)}
\end{equation}

Therefore

\begin{equation}\label{RP3}
 \langle {\bf p},{\bf p} \rangle
\begin{cases}
>0 & \text{when \;$0 < \psi<\dfrac{2}{3} \pi$},\\
=0 & \text{when \;$\psi= \dfrac{2}{3} \pi$},\\
<0 & \text{when \;$\dfrac{2}{3}\pi< \psi < \omega$}.
\end{cases}
\end{equation}

\begin{figure}[H]
\includegraphics[height=70mm,width=90mm]{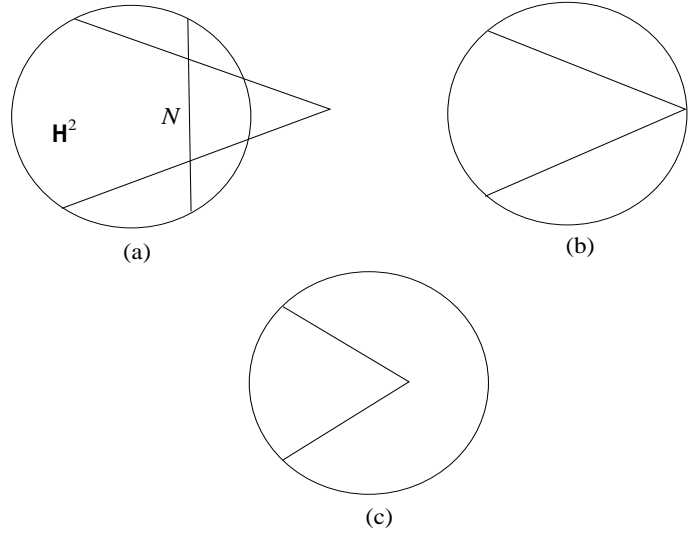}
\centering
\caption{(a) $\langle {\bf p},{\bf p}\rangle>0$ : lines meet outside
  $\Bbb H^2$. (b) $\langle {\bf p},{\bf p}\rangle=0$ : lines meet on
  $\partial\Bbb H^2$. (c) $\langle {\bf p},{\bf p}\rangle<0$ : lines meet
  inside $\Bbb H^2$.}\label{fig:RP2}
\end{figure}

Information contained in (\ref{RP3}) corresponds to the configurations
depicted in figure \ref{fig:RP2} in the case of $\Bbb H^2$ embedded in
$\Bbb RP^2$ and is true for $n$ hyperplanes in $\Bbb H^n$ embedded in $\Bbb RP^n$.

Referring to (\ref{N}), figure \ref{fig:fig6} and figure \ref{fig:RP2}, when $\langle {\bf p},{\bf
  p}\rangle>0$ there is a plane $N$ inside $\Bbb
H^3$. When $\langle {\bf p},{\bf p}\rangle=0$ the plane $N$ is the point
of intersection of the planes $A, C$ and $D$ on $\partial\Bbb
H^3$. When $\langle {\bf p},{\bf p}\rangle<0$ the plane $N$ is the
point of intersection of the planes $A, C$ and $D$ inside $\Bbb H^3$.

We now verify that $P_\psi \subset \Bbb H^3 \cup \partial\Bbb H^3$
for $\psi\in[0,\frac{2}{3}\pi]$. Since $P_\psi$ is symmetric with
respect to the origin it is sufficient to show that the vertices $(A
\cap N \cap C)$ and $(A \cap N \cap D)$ of figure \ref{fig:fig6} live
inside $\Bbb H^3 \cup \partial\Bbb H^3$.
Let 
\begin{equation}\label{ANC}
f(\psi)=R^2-(A \cap N \cap C)^2 = \frac{1}{2}(R^2-c^2)(2+c^2-R^2)
\end{equation}
where
$(A \cap N \cap C)^2$ denotes square of the distance of vertex $(A\cap N
\cap C)$ from the origin.
We have $f>0$ for $\psi \in (0,\frac{2}{3}\pi)$ and $f=0$ when
$\psi=0$ or $\frac{2}{3}\pi$.

Let
\begin{equation}\label{AND}
g(\psi)=R^2-(A\cap N \cap D)^2= \frac{(R^2-1)(2+c^2-R^2)}{c^2+1}
\end{equation}

We note $g>0$ when $\psi \in (0,\frac{2}{3}\pi)$ and $g=0$ if $\psi=0$
or $\frac{2}{3}\pi$.

Therefore $P_\psi \subset {\Bbb H^3 \cup \partial
  \Bbb H^3}$ when $\psi \in [0,\frac{2}{3}\pi]$. Hence $M_\psi$ is
  hyperbolic when $\psi \in [0,\frac{2}{3}\pi]$. 

\begin{figure}[H]
\includegraphics[height=50mm,width=50mm]{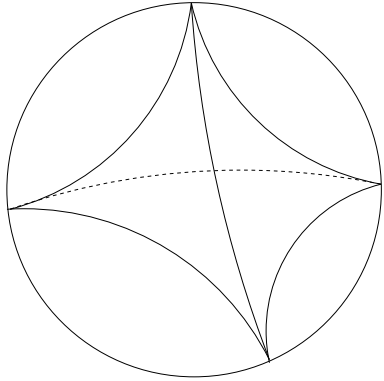}
\centering
\caption{$P_{\frac{2}{3}\pi}$ inside $\Bbb H^3 \cup \partial \Bbb H^3$}\label{fig:P23Pi}
\end{figure}

\begin{figure}[H]
\includegraphics[height=40mm,width=60mm]{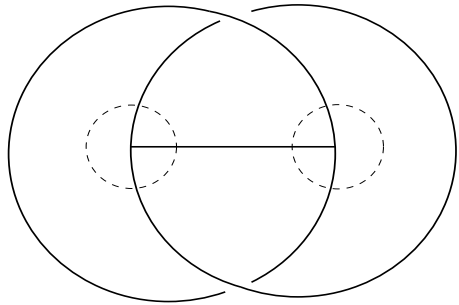}
\centering
\caption{Singular set of $M_\psi$ when
  $\psi\in[\frac{2}{3}\pi,\omega)$ and its cusp neighbourhoods when $\psi=\frac{2}{3}\pi$.}\label{fig:S23Pi}
\end{figure}

\begin{figure}[H]
\includegraphics[height=27mm,width=50mm]{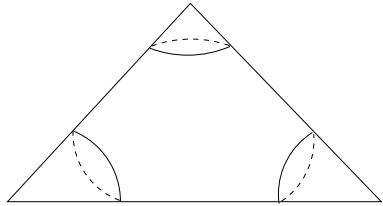}
\centering
\caption{A Euclidean turnover}\label{fig:Eturnover}
\end{figure}

We note from (\ref{RP3}) and figure \ref{fig:RP2} that $M_{\frac{2}{3}\pi}$
  is cusped. From figure \ref{fig:P23Pi} we
  observe that $M_{\frac{2}{3}\pi}$ has two cusps, it is topologically
  $S^3$ with singular set as shown in figure \ref{fig:S23Pi} and its
  cusp neighborhoods are Euclidean turnovers as in figure \ref{fig:Eturnover}.

From (\ref{RP3}) and figure \ref{fig:RP2} we observe that the bounding
planes of $P_\psi$ meet inside $\Bbb H^3$ when
$\psi\in(\frac{2}{3}\pi,\omega)$. Therefore, $P_\psi$ is in $\Bbb
H^3$. Hence $M_\psi$ is hyperbolic when $\psi\in(\frac{2}{3}\pi,\omega)$. The singular set of $M_\psi$ for $\psi\in(\frac{2}{3}\pi,\omega)$ is shown in figure \ref{fig:S23Pi}.

\begin{figure}[H]
\includegraphics[height=40mm,width=70mm]{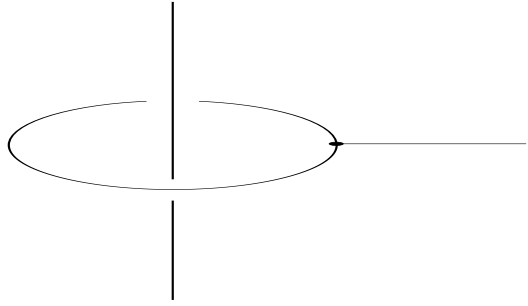}
\centering
\caption{Singular set of $M_\omega$}\label{fig:Somega}
\end{figure}

We note from (\ref{R2}) that $R(\omega)=\infty$ hence $M_\omega$ is
Euclidean. $P_\omega$ is a tetrahedron with deleted vertices
from which $M_\omega$ is obtained using identifications. We observe that $M_\omega$ is $R^3$. The singular set of $M_\omega$ is shown in figure \ref{fig:Somega}. 

\section{Cone angles larger than $\omega$}

\begin{prop}
$M_\psi$ is spherical when cone angle $\psi$ is larger than
$\omega$. There exists $\zeta\in \Bbb R$ such that $M_\psi$ is
topologically $S^3$ with a self intersecting singular set when $\omega
< \psi <\zeta$. $M_\psi$ is the suspension of a sphere with four cone points when $\zeta \leq \psi$.
\end{prop} 

{\bf Proof} \;When $\omega<\psi<\pi$ we note $R^2$ the square radius of the
Klein model becomes negative so that the model has imaginary
radius. This leads us to use the spherical bilinear form for ${\bf
  v}=(v_1,v_2,v_3,v_4),{\bf w}=(w_1,w_2,w_3,w_4)\in \Bbb R^4$ :

\begin{equation}\label{SBF}
{\langle {\bf v},{\bf w} \rangle}_\Bbb S = v_1w_1+v_2w_2+v_3w_3+v_4w_4
\end{equation}

when $\psi\in(\omega,\pi)$.

This is the usual scalar product on $\Bbb R^4$. Let ${\bf
  x}=(x_1,x_2,x_3,x_4)\in \Bbb R^4$. Since ${\Bbb H}^3_R=\{x \mid {\langle
  {\bf x},{\bf x} \rangle}_\Bbb H = -R^2\}$ we define ${\Bbb S}^3_R=\{x \mid {\langle
  {\bf x},{\bf x} \rangle}_\Bbb S = R^2\}$. Thus ${\Bbb S}^3_R$ is the
  sphere of radius $R$.
The ``Klein model'' ${\Bbb K}^3_R$ for the sphere of radius $R$ is the
  hyperplane ${\Bbb K}^3_R = \{{\bf x} \mid t=R\}$. In contrast to the
  hyperbolic case we don't need to verify that the polytope lies
  inside the sphere of radius $R$ in ${\Bbb K}^3_R$ since projection
  from the origin which is not conformal defines a 1-1 correspondence
  between ${\Bbb K}^3_R$ and the upper hemisphere of ${\Bbb
  S}^3_R$. The metric on ${\Bbb K}^3_R$ is then the pull back of the
  metric on ${\Bbb S}^3_R$. As in the hyperbolic case reflections in
  planes through the origin and rotation about axes through the origin
  are both Euclidean and spherical isometries.

If a plane in ${\Bbb K}^3_R$ has equation $\alpha x+\beta y+\gamma
z=\delta$ and $t=R$ then in homogeneous coordinates the equation is
$\alpha x+\beta y+\gamma z-(\delta /R)t=0$. In the spherical case the
pole has homogeneous coordinates $(\alpha,\beta,\gamma,-\delta
/R)$ whereas in the hyperbolic case the pole has coordinates $(\alpha,\beta,\gamma,\delta /R)$.
 If a pair of planes with poles ${\bf v}$ and ${\bf w}$ in the
 spherical case intersect with dihedral angle $\theta$ then
\begin{equation}\label{Sangle}
\pm\cos\theta=\dfrac{{\langle{\bf v},{\bf w}\rangle}_\Bbb S}{\sqrt{{\langle{\bf
      v},{\bf v}\rangle}_\Bbb S {\langle{\bf w},{\bf w}\rangle}_\Bbb S}}
\end{equation}

Using (\ref{SBF}), (\ref{Sangle}) and the spherical version of (\ref{planes}) we obtain
\begin{equation}
0\leq R^2 = -\dfrac{1+\cos \psi}{2\cos \dfrac{\psi}{4}-1+\cos \psi}
\end{equation}

when $\omega\leq\psi\leq \zeta$.

\begin{figure}[H]
\includegraphics[height=59mm,width=63mm]{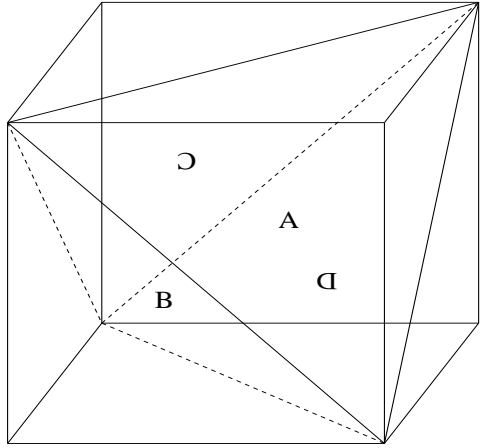}
\centering
\caption{$P_\psi$ inside the reference box}\label{fig:tetra}
\end{figure}

When $\psi \in (\omega,\pi)$ we observe that $P_\psi$ is a tetrahedron
with identification as in figure \ref{fig:tetra}. 

Therefore $M_\psi$ is $S^3$ with a self intersecting singular set as in figure
\ref{fig:S23Pi} with spherical turnover neighbourhoods.

\begin{figure}[H]
\includegraphics[height=35mm,width=40mm]{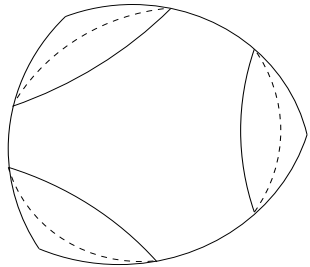}
\centering
\caption{A spherical turnover}\label{fig:turnover}
\end{figure}

\begin{figure}[H]
\includegraphics[height=35mm,width=85mm]{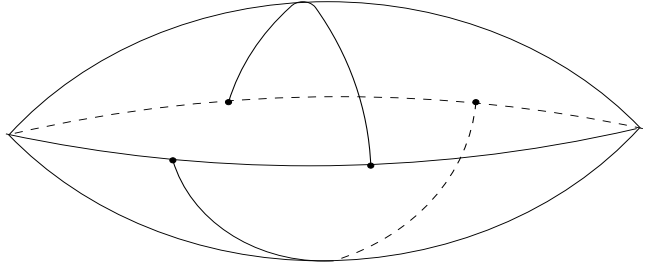}
\centering
\caption{$P_\pi$ is a lens with angle $\dfrac{\pi}{4}$}\label{fig:PPI}
\end{figure}

Since $R^2(\pi)=0=c^2(\pi)$ the three dimensional model has collapsed
into a two dimensional disc at cone angle $\pi$. By projection onto
$S_1^3=\{{\bf x} \mid x^2+y^2+z^2+t^2=1\}$ the sphere of radius 1 we
observe that $P_\pi$ is a lens with angle $\psi/4$ as in figure \ref{fig:PPI} from which $M_\pi$ is obtained by identifications.

\begin{figure}[H]
\includegraphics[height=40mm,width=70mm]{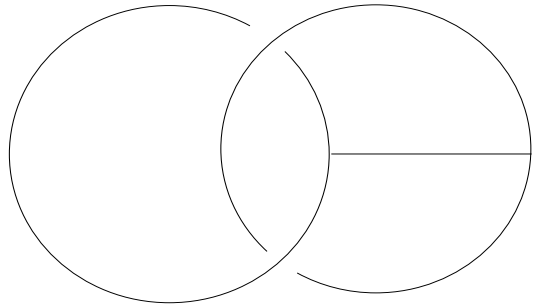}
\centering
\caption{Singular set of $M_\pi$}\label{fig:SgPi}
\end{figure}

We note $\zeta=4\cos^{-1}(\dfrac{2}{\sqrt3}\cos(\dfrac{1}{3}\cos^{-1}(\dfrac{-3\sqrt3}{8}+\dfrac{4}{3}\pi)))\approx
5.191298...$ is the second smallest value of $\psi$ for which $R$ is
infinite. The equation of $c^2$ is the same as in the hyperbolic case.
The model continues to be $S^3$ for $\pi<\psi<\zeta$ therefore
$M_\psi$ is spherical with a self intersecting singular set as in
figure \ref{fig:SgPi} when $\psi \in (\pi,\zeta)$.

\begin{figure}[H]
\includegraphics[height=40mm,width=80mm]{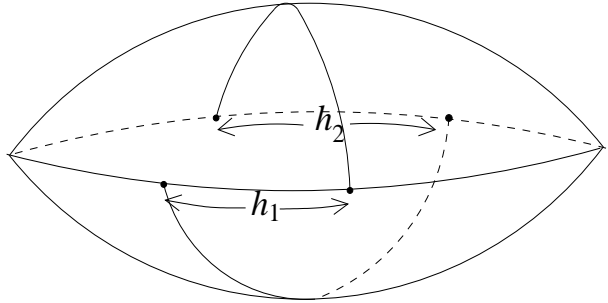}
\centering
\caption{Segments $h_1$ and $h_2$ decrease in length as $\psi$
  increases towards $\zeta$}\label{fig:Ppih}
\end{figure}

\begin{figure}[H]
\includegraphics[height=40mm,width=80mm]{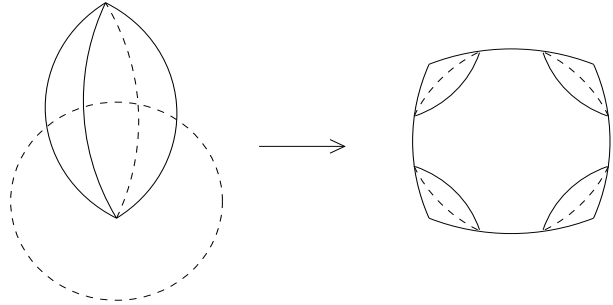}
\centering
\caption{$M_\zeta$ is the suspension of a sphere with four cone points}\label{fig:Pzeta}
\end{figure}

As cone angle $\psi$ approaches $\zeta$ the segments $h_1$ and $h_2$
decrease in length towards $0$ as shown in figure \ref{fig:Ppih} so that we get
the singular set shown in figure \ref{fig:Pzeta} when $\psi=\zeta$. 

We therefore obtain $M_\psi$ as the suspension of a sphere with four
cone points an in figure \ref{fig:Pzeta} when $\zeta<\psi$.

\section{Spontaneous surgery}

\begin{figure}[H]
\includegraphics[height=60mm,width=80mm]{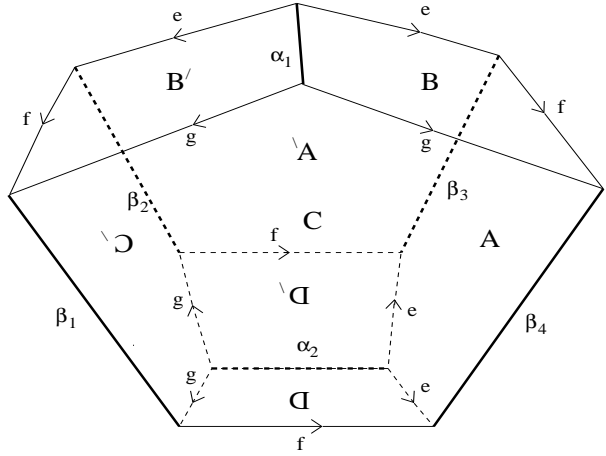}
\centering
\caption{$P_\psi$ after spontaneous surgery gives $S^3$ with the
  Whitehead link as its singular set}\label{fig:wlpoly}
\end{figure}

We now look at how we get spontaneous surgery and the Whitehead link
as the singular set.
$P_\psi$ and its identifications ``after''spontaneous surgery are
shown in figure \ref{fig:wlpoly}. Face pairings are
$A \longleftrightarrow A', B\longleftrightarrow B',
C\longleftrightarrow C'$ and $D \longleftrightarrow D'$ and edges
with the same label are identified.

We note that $M_\psi$ ``after'' spontaneous surgery is the result of
performing $(0,1)$-Dehn surgery on one component of the Whitehead link
in $S^3$. Dihedral angles at the $\beta$ edges of $P_\psi$ are now the cone angle
$\psi$ and other incidence angles between the planes of $P_\psi$ are
as they were prior to spontaneous surgery.
Calculations to show the behaviour of $M_\psi$ after surgery remain to
be done.


\input{wl1.bbl}

\;
\;
\;
Dr Aalam\\
PO Box 18810\\
London SW7 2ZR\\
UK.\\
email: aalam@mth.kcl.ac.uk

\end{document}

%% file: wl1.bbl
\providecommand{\bysame}{\leavevmode\hbox to3em{\hrulefill}\thinspace}
\providecommand{\MR}{\relax\ifhmode\unskip\space\fi MR }
\providecommand{\MRhref}[2]{%
  \href{http://www.ams.org/mathscinet-getitem?mr=#1}{#2}
}
\providecommand{\href}[2]{#2}

%% file: wl1.bbl
\begin{thebibliography}{1}

\bibitem{HLM}
H.~Hilden{,}~M. Lozano and J.~Montesinos, \emph{On a {R}emarkable {P}olyhedron
  {G}eometrising the {F}igure {E}ight {K}not {C}one {M}anifolds}, J. Math. Sci.
  Univ. Tokyo \textbf{2} (1995), no.~3, 501--561.

\bibitem{Rat}
J.~G. Ratcliffe, \emph{Foundations of {H}yperbolic {M}anifolds},
  Springer-Verlag, 1994.

\bibitem{Thu2}
W.~P. Thurston, \emph{The {G}eometry and {T}opology of {T}hree-{M}anifolds},
  Princeton Univ. Math. Dept., 1978, available from
  http://msri.org/publications/books/gt3m/.

\bibitem{Thu1}
\bysame, \emph{Three-{D}imensional {G}eometry and {T}opology}, vol.~1,
  Princeton Univ. Press, 1997.

\end{thebibliography}
